# STEIN'S METHOD AND EXACT BERRY–ESSEEN ASYMPTOTICS FOR FUNCTIONALS OF GAUSSIAN FIELDS

By Ivan Nourdin and Giovanni Peccati

*Université Paris VI and Université Paris Ouest*

We show how to detect optimal Berry–Esseen bounds in the normal approximation of functionals of Gaussian fields. Our techniques are based on a combination of Malliavin calculus, Stein's method and the method of moments and cumulants, and provide de facto local (one-term) Edgeworth expansions. The findings of the present paper represent a further refinement of the main results proven in Nourdin and Peccati [*Probab. Theory Related Fields* **145** (2009) 75–118]. Among several examples, we discuss three crucial applications: (i) to Toeplitz quadratic functionals of continuous-time stationary processes (extending results by Ginovyan [*Probab. Theory Related Fields* **100** (1994) 395–406] and Ginovyan and Sahakyan [*Probab. Theory Related Fields* **138** (2007) 551–579]); (ii) to "exploding" quadratic functionals of a Brownian sheet; and (iii) to a continuous-time version of the Breuer–Major CLT for functionals of a fractional Brownian motion.

**1. Introduction.** Let $\{F_n : n \geq 1\}$ be a sequence of zero-mean real-valued random variables and consider a standard Gaussian variable $N \sim \mathcal{N}(0,1)$. Assume that each $F_n$ is a functional of an infinite-dimensional Gaussian field and suppose that, as $n \to \infty$,

$$F_n \xrightarrow{\text{Law}} N. \tag{1.1}$$

In the paper [20], the present authors demonstrated that one can naturally combine Malliavin calculus (see, e.g., [13, 21]) with Stein's method (see, e.g., [4, 29, 33, 34]) in order to obtain explicit bounds of the type

$$d(F_n, N) \leq \varphi(n), \qquad n \geq 1, \tag{1.2}$$

where $d(F_n, N)$ stands for some appropriate distance (e.g., the Kolmogorov distance or the total variation distance) between the laws of $F_n$ and $N$, and $\varphi(n)$ is some positive sequence converging to zero. The aim of the present work is to develop several techniques, allowing us to assess the optimality of the bound $\varphi(n)$ appearing in (1.2) for a given sequence $\{F_n\}$. Formally, one says that the bound $\varphi(n)$ is *optimal* for the sequence $\{F_n\}$ and the distance $d$ whenever there exists a constant $c \in (0,1)$ (independent of $n$) such that, for $n$ sufficiently large,

$$(1.3) \qquad c < d(F_n, N)/\varphi(n) \leq 1.$$

We shall establish relations such as (1.3) by pushing the Malliavin-type approach to Stein's method (initiated in [20]) one step further. In particular, the findings of this paper represent a new and substantial refinement of the central limit theorems (CLTs) for functionals of Gaussian fields which were proven in [22, 23, 25, 26]. Once again, our techniques do not require that the random variables $\{F_n\}$ have the specific form of partial sums. Indeed, we will see, in Sections 4–6 below, that our results yield optimal Berry–Esseen-type bounds for CLTs involving objects as diverse as: (i) Toeplitz quadratic functionals of continuous-time stationary processes; (ii) quadratic functionals of a Brownian motion or of a Brownian sheet indexed by a compact set of $\mathbb{R}^d$ ($d \geq 2$); and (iii) polynomial functionals constructed from a fractional Brownian motion.

Note that, in the subsequent sections, we shall focus uniquely on the normal approximation of random variables with respect to the *Kolmogorov distance*. This distance is defined as

$$(1.4) \qquad d_{\mathrm{Kol}}(X, Y) = \sup_{z \in \mathbb{R}} |P(X \leq z) - P(Y \leq z)|$$

for any pair of random variables $X$ and $Y$. It will later become clear that many results of the present paper extend almost verbatim to alternate distances, such as the Wasserstein and the total variation distance, between laws of real-valued random variables.

Our basic approach can be described as follows. Fix $z \in \mathbb{R}$ and consider the *Stein equation*

$$(1.5) \qquad \mathbf{1}_{(-\infty, z]}(x) - \Phi(z) = f'(x) - x f(x), \qquad x \in \mathbb{R},$$

where, here, and for the rest of the paper, we use the standard notation $\Phi(z) = P(N \leq z)$ [$N \sim \mathcal{N}(0,1)$] and where $\mathbf{1}_A$ stands for the indicator of a set $A$. It is well known that, for every fixed $z$, (1.5) admits a solution $f_z$ such that $\|f_z\|_\infty \leq \sqrt{2\pi}/4$ and $\|f_z'\|_\infty \leq 1$ (see, e.g., [4], Lemma 2.2, or formulae (2.10)–(2.11) below). Now, suppose that the elements of the sequence $\{F_n\}$ appearing in (1.1) are functionals of some Gaussian field, say $X$, and assume that each $F_n$ is differentiable in the sense of Malliavin calculus (see Section 2.1 for details). Denote by $DF_n$ the Malliavin derivative of $F_n$ and by $L^{-1}$ the *pseudo-inverse of the Ornstein–Uhlenbeck generator* (again, see Section 2.1). Recall that $DF_n$ is a random element with values in an appropriate Hilbert space $\mathfrak{H}$. In [20], Section 3, we proved and applied the following relations, direct consequences of the fact that $f_z$ solves (1.5) on the one hand, and of the celebrated integration by parts formula of Malliavin calculus on the other hand: for every $z \in \mathbb{R}$,

$$
\begin{aligned}
P(F_n \leq z) - \Phi(z) &= E[f'_z(F_n) - F_n f_z(F_n)] \\
&= E[f'_z(F_n)(1 - \langle DF_n, -DL^{-1}F_n \rangle_{\mathfrak{H}})].
\end{aligned}
\tag{1.6}
$$

By using (1.4), applying the Cauchy–Schwarz inequality to the right-hand side of (1.6) and using the fact that $f'_z$ is bounded by 1, one immediately obtains that

$$
d_{\mathrm{Kol}}(F_n, N) \leq \sqrt{E[(1 - \langle DF_n, -DL^{-1}F_n \rangle_{\mathfrak{H}})^2]}.
\tag{1.7}
$$

The starting point of [20] was that, in several crucial cases (e.g., when each $F_n$ is a multiple Wiener–Itô integral of a fixed order), the upper bound

$$
\varphi(n) := \sqrt{E[(1 - \langle DF_n, -DL^{-1}F_n \rangle_{\mathfrak{H}})^2]}, \qquad n \geq 1,
\tag{1.8}
$$

is such that: (i) the quantity $\varphi(n)$ can be explicitly computed (e.g., in terms of contraction operators); (ii) $\varphi(n) \to 0$ as $n \to \infty$; and (iii) $\varphi(n)$ is directly related to quantities playing a fundamental role in the CLTs for functionals of Gaussian fields proven in [22, 23, 25, 26]. The aim of the present paper is to establish conditions on the sequence $\{F_n\}$ ensuring that the ratios

$$
\frac{E[f'_z(F_n)(1 - \langle DF_n, -DL^{-1}F_n \rangle_{\mathfrak{H}})]}{\varphi(n)}, \qquad n \geq 1,
\tag{1.9}
$$

involving (1.8) and the right-hand side of (1.6), converge to a nonzero limit for all $z$ outside some finite set. Such a result immediately yields the existence of a constant $c$, verifying (1.3) for $d = d_{\mathrm{Kol}}$. We will show that a very effective way to prove the convergence of the quantities appearing in (1.9) is to characterize the joint convergence in distribution of the random vectors

$$
\left( F_n, \frac{1 - \langle DF_n, -DL^{-1}F_n \rangle_{\mathfrak{H}}}{\varphi(n)} \right), \qquad n \geq 1,
\tag{1.10}
$$

toward a two-dimensional Gaussian vector with nonzero covariance. The applications presented in Sections 4–6 will show that this specific convergence takes place in several crucial situations, involving, for instance, quadratic or polynomial functionals of stationary Gaussian processes. We will see that, in order to prove a CLT for the vector appearing in (1.10), a useful tool is the multidimensional version of the CLT for multiple stochastic integrals which was proven in [26]. Also, it is interesting to note that if each $F_n$ in (1.1) is a double stochastic integral, then our conditions can be expressed exclusively in terms of the second, third, fourth and eighth cumulants associated with the sequence $\{F_n\}$; see Section 3.3 below.

The rest of the paper is organized as follows. Section 2 deals with preliminaries concerning Malliavin calculus, Stein's method and related topics. Section 3 contains our main results, with special attention devoted to random variables belonging to the second Wiener chaos of a Gaussian field. In Section 4, we develop an application to Toeplitz quadratic functionals of stationary continuous-time Gaussian processes, thus extending and refining some results by Ginovyan [7] and Ginovyan and Sahakyan [8]. Section 5 is devoted to quadratic functionals of Brownian motion and of the Brownian sheet, whereas Section 6 focuses on a continuous-time version of the Breuer–Major CLT for processes subordinated to a fractional Brownian motion.

## 2. Preliminaries.

2.1. *Gaussian fields and Malliavin calculus.* We shall now provide a short description of the tools of Malliavin calculus that will be needed in the forthcoming sections. The reader is referred to the monographs [13] and [21] for any unexplained concepts or results.

Let $\mathfrak{H}$ be a real separable Hilbert space. We denote by $X = \{X(h) : h \in \mathfrak{H}\}$ an isonormal Gaussian process over $\mathfrak{H}$. By definition, $X$ is a centered Gaussian family indexed by the elements of $\mathfrak{H}$ and such that, for every $h, g \in \mathfrak{H}$,

$$(2.1) \qquad E[X(h)X(g)] = \langle h, g \rangle_{\mathfrak{H}}.$$

In what follows, we shall use the notation $L^2(X) = L^2(\Omega, \sigma(X), P)$. For every $q \geq 1$, we write $\mathfrak{H}^{\otimes q}$ to indicate the $q$th tensor power of $\mathfrak{H}$; the symbol $\mathfrak{H}^{\odot q}$ stands for the $q$th *symmetric* tensor power of $\mathfrak{H}$, equipped with the norm $\sqrt{q!}\|\cdot\|_{\mathfrak{H}^{\otimes q}}$. We denote by $I_q$ the isometry between $\mathfrak{H}^{\odot q}$ and the $q$th Wiener chaos of $X$. It is well known (again, see [21], Chapter 1, or [13]) that any random variable $F$ belonging to $L^2(X)$ admits the chaotic expansion

$$(2.2) \qquad F = \sum_{q=0}^{\infty} I_q(f_q),$$

where $I_0(f_0) := E[F]$, the series converges in $L^2$ and the kernels $f_q \in \mathfrak{H}^{\odot q}$, $q \geq 1$, are uniquely determined by $F$. In the particular case where $\mathfrak{H} = L^2(A, \mathscr{A}, \mu)$, where $(A, \mathscr{A})$ is a measurable space and $\mu$ is a $\sigma$-finite and nonatomic measure, one has that $\mathfrak{H}^{\odot q} = L_s^2(A^q, \mathscr{A}^{\otimes q}, \mu^{\otimes q})$ is the space of symmetric and square-integrable functions on $A^q$. Moreover, for every $f \in \mathfrak{H}^{\odot q}$, $I_q(f)$ coincides with the multiple Wiener–Itô integral (of order $q$) of $f$ with respect to $X$ (see [21], Chapter 1). It is well known that a random variable of the type $I_q(f)$, $f \in \mathfrak{H}^{\odot q}$, has finite moments of all orders (see, e.g., [13], Chapter VI). Moreover, any nonzero finite sum of multiple stochastic integrals has a law which is absolutely continuous with respect to Lebesgue measure (see, e.g., Shigekawa [32] for a proof of this fact; see [21], Chapter 1, or [30] for a connection between multiple Wiener–Itô integrals and Hermite polynomials on the real line). For every $q \geq 0$, we denote by $J_q$ the orthogonal projection operator on the $q$th Wiener chaos associated with $X$ so that, if $F \in L^2(\sigma(X))$ is as in (2.2), then $J_q F = I_q(f_q)$ for every $q \geq 0$.

Let $\{e_k, k \geq 1\}$ be a complete orthonormal system in $\mathfrak{H}$. Given $f \in \mathfrak{H}^{\odot p}$ and $g \in \mathfrak{H}^{\odot q}$, for every $r = 0, \ldots, p \wedge q$, the $r$th contraction of $f$ and $g$ is the element of $\mathfrak{H}^{\otimes(p+q-2r)}$ defined as

$$(2.3) \quad f \otimes_r g = \sum_{i_1, \ldots, i_r = 1}^{\infty} \langle f, e_{i_1} \otimes \cdots \otimes e_{i_r} \rangle_{\mathfrak{H}^{\otimes r}} \otimes \langle g, e_{i_1} \otimes \cdots \otimes e_{i_r} \rangle_{\mathfrak{H}^{\otimes r}}.$$

In the particular case where $\mathfrak{H} = L^2(A, \mathscr{A}, \mu)$ (with $\mu$ nonatomic), one has that

$$f \otimes_r g = \int_{A^r} f(t_1, \ldots, t_{p-r}, s_1, \ldots, s_r)$$
$$\times g(t_{p-r+1}, \ldots, t_{p+q-2r}, s_1, \ldots, s_r) \, d\mu(s_1) \cdots d\mu(s_r).$$

Moreover, $f \otimes_0 g = f \otimes g$ equals the tensor product of $f$ and $g$ while, for $p = q$, $f \otimes_p g = \langle f, g \rangle_{\mathfrak{H}^{\otimes p}}$. Note that, in general (and except for trivial cases), the contraction $f \otimes_r g$ is *not* a symmetric element of $\mathfrak{H}^{\otimes(p+q-2r)}$. The canonical symmetrization of $f \otimes_r g$ is written $f \widetilde{\otimes}_r g$. We also have the following multiplication formula: if $f \in \mathfrak{H}^{\odot p}$ and $g \in \mathfrak{H}^{\odot q}$, then

$$(2.4) \qquad I_p(f) I_q(g) = \sum_{r=0}^{p \wedge q} r! \binom{p}{r} \binom{q}{r} I_{p+q-2r}(f \widetilde{\otimes}_r g).$$

Let $\mathscr{S}$ be the set of all smooth cylindrical random variables of the form

$$F = g(X(\phi_1), \ldots, X(\phi_n)),$$

where $n \geq 1$, $g : \mathbb{R}^n \to \mathbb{R}$ is a smooth function with compact support and $\phi_i \in \mathfrak{H}$. The Malliavin derivative of $F$ with respect to $X$ is the element of $L^2(\Omega, \mathfrak{H})$ defined as

$$DF = \sum_{i=1}^{n} \frac{\partial g}{\partial x_i}(X(\phi_1), \ldots, X(\phi_n)) \phi_i.$$

Also, $DX(h) = h$ for every $h \in \mathfrak{H}$. By iteration, one can define the $m$th derivative $D^m F$ [which is an element of $L^2(\Omega, \mathfrak{H}^{\otimes m})$] for every $m \geq 2$. As usual, for $m \geq 1$, $\mathbb{D}^{m,2}$ denotes the closure of $\mathscr{S}$ with respect to the norm $\|\cdot\|_{m,2}$, defined by the relation

$$\|F\|_{m,2}^2 = E[F^2] + \sum_{i=1}^{m} E[\|D^i F\|_{\mathfrak{H}^{\otimes i}}^2].$$

Note that if $F$ is equal to a finite sum of multiple Wiener–Itô integrals, then $F \in \mathbb{D}^{m,2}$ for every $m \geq 1$. The Malliavin derivative $D$ verifies the following chain rule: if $\varphi : \mathbb{R}^n \to \mathbb{R}$ is in $\mathscr{C}_b^1$ (i.e., the collection of bounded continuously differentiable functions with a bounded derivative) and if $\{F_i\}_{i=1,\ldots,n}$ is a vector of elements of $\mathbb{D}^{1,2}$, then $\varphi(F_1,\ldots,F_n) \in \mathbb{D}^{1,2}$ and

$$D\varphi(F_1,\ldots,F_n) = \sum_{i=1}^{n} \frac{\partial \varphi}{\partial x_i}(F_1,\ldots,F_n) DF_i.$$

Observe that the previous formula still holds when $\varphi$ is a Lipschitz function and the law of $(F_1,\ldots,F_n)$ has a density with respect to the Lebesgue measure on $\mathbb{R}^n$ (see, e.g., Proposition 1.2.3 in [21]). We denote by $\delta$ the adjoint of the operator $D$, also called the *divergence operator*. A random element $u \in L^2(\Omega, \mathfrak{H})$ belongs to the domain of $\delta$, noted $\mathrm{Dom}\,\delta$, if and only if it verifies

$$|E\langle DF, u\rangle_\mathfrak{H}| \leq c_u \|F\|_{L^2} \qquad \text{for any } F \in \mathscr{S},$$

where $c_u$ is a constant depending uniquely on $u$. If $u \in \mathrm{Dom}\,\delta$, then the random variable $\delta(u)$ is defined by the duality relationship (i.e., the "integration by parts formula")

(2.5) $$E(F\delta(u)) = E\langle DF, u\rangle_\mathfrak{H},$$

which holds for every $F \in \mathbb{D}^{1,2}$.

The operator $L$, acting on square-integrable random variables of the type (2.2), is defined through the projection operators $\{J_q\}_{q\geq 0}$ as $L = \sum_{q=0}^{\infty} -q J_q$ and is called the *infinitesimal generator of the Ornstein–Uhlenbeck semigroup*. It verifies the following crucial property: a random variable $F$ is an element of $\mathrm{Dom}\,L$ $(= \mathbb{D}^{2,2})$ if and only if $F \in \mathrm{Dom}\,\delta D$ (i.e., $F \in \mathbb{D}^{1,2}$ and $DF \in \mathrm{Dom}\,\delta$) and, in this case, $\delta DF = -LF$. Note that a random variable $F$ as in (2.2) is in $\mathbb{D}^{1,2}$ (resp., $\mathbb{D}^{2,2}$) if and only if

$$\sum_{q=1}^{\infty} q\|f_q\|_{\mathfrak{H}^{\odot q}}^2 < \infty \qquad \left(\text{resp.}, \sum_{q=1}^{\infty} q^2 \|f_q\|_{\mathfrak{H}^{\odot q}}^2 < \infty\right),$$

and also $E[\|DF\|_\mathfrak{H}^2] = \sum_{q\geq 1} q\|f_q\|_{\mathfrak{H}^{\odot q}}^2$. If $\mathfrak{H} = L^2(A, \mathscr{A}, \mu)$ (with $\mu$ nonatomic), then the derivative of a random variable $F$ as in (2.2) can be identified with the element of $L^2(A \times \Omega)$ given by

(2.6) $$D_a F = \sum_{q=1}^{\infty} q I_{q-1}(f_q(\cdot, a)), \qquad a \in A.$$

We also define the operator $L^{-1}$, which is the *pseudo-inverse* of $L$, as follows: for every $F \in L^2(X)$, we set $L^{-1}F = \sum_{q \geq 1} \frac{1}{q} J_q(F)$. Note that $L^{-1}$ is an operator with values in $\mathbb{D}^{2,2}$ and that $LL^{-1}F = F - E(F)$ for all $F \in L^2(X)$.

The following lemma generalizes Lemma 2.1 in [19].

LEMMA 2.1. *Let $F \in \mathbb{D}^{1,2}$ be such that $E(F) = 0$. Suppose that, for some integer $s \geq 0$, $E|F|^{s+2} < \infty$. Then,*

$$(2.7) \qquad E(F^s \langle DF, -DL^{-1}F \rangle_{\mathfrak{H}}) = \frac{1}{s+1} E(F^{s+2}).$$

PROOF. Since $L^{-1}F \in \mathbb{D}^{2,2}$, we can write

$$E(F^s \langle DF, -DL^{-1}F \rangle_{\mathfrak{H}})$$
$$= \frac{1}{s+1} E(\langle D(F^{s+1}), D(-L^{-1}F) \rangle_{\mathfrak{H}})$$
$$= -\frac{1}{s+1} E(\delta DL^{-1}F \times F^{s+1}) \qquad \text{[by integration by parts (2.5)]}$$
$$= \frac{1}{s+1} E(F^{s+2}) \qquad \text{(by the relation } -\delta DL^{-1}F = F\text{)}. \qquad \square$$

REMARK 2.2. If $F = I_q(f)$ for some $q \geq 2$ and $f \in \mathfrak{H}^{\odot q}$, then

$$(2.8) \quad \langle DF, -DL^{-1}F \rangle_{\mathfrak{H}} = \langle DI_q(f), -DL^{-1}I_q(f) \rangle_{\mathfrak{H}} = \frac{1}{q} \|DI_q(f)\|_{\mathfrak{H}}^2$$

so that (2.7) yields, for every integer $s \geq 1$, that

$$(2.9) \qquad E(I_q(f)^s \|DI_q(f)\|_{\mathfrak{H}}^2) = \frac{q}{s+1} E(I_q(f)^{s+2}).$$

2.2. *Stein's method and normal approximation on a Gaussian space.* We start by recalling that, for every fixed $z \in \mathbb{R}$, the function

$$(2.10) \qquad f_z(x) = e^{x^2/2} \int_{-\infty}^{x} [\mathbf{1}_{(-\infty,z]}(a) - \Phi(z)] e^{-a^2/2} \, da$$

$$(2.11) \qquad = \begin{cases} \sqrt{2\pi} e^{x^2/2} \Phi(x)(1 - \Phi(z)), & \text{if } x \leq z, \\ \sqrt{2\pi} e^{x^2/2} \Phi(z)(1 - \Phi(x)), & \text{if } x > z, \end{cases}$$

is a solution to the Stein equation (1.5), also verifying $\|f_z\|_\infty \leq \sqrt{2\pi}/4$ and $\|f_z'\|_\infty \leq 1$.

The following lemma will play a crucial role in the sequel; see also (1.6). Its content is the starting point of [20].

LEMMA 2.3. *Let $F \in \mathbb{D}^{1,2}$ have zero mean. Assume, moreover, that $F$ has an absolutely continuous law with respect to the Lebesgue measure. Then, for every $z \in \mathbb{R}$,*

$$P(F \leq z) - \Phi(z) = E[f'_z(F)(1 - \langle DF, -DL^{-1}F \rangle_{\mathfrak{H}})].$$

PROOF. Fix $z \in \mathbb{R}$. Since $f_z$ solves the Stein equation (1.5), we have $P(F \leq z) - \Phi(z) = E[f'_z(F) - Ff_z(F)]$. Now, observe that one can write $F = LL^{-1}F = -\delta DL^{-1}F$. By using the integration by parts formula (2.5) and the fact that $Df_z(F) = f'_z(F)DF$ (note that, for this formula to hold with $f_z$ only Lipschitz, one needs $F$ to have an absolutely continuous law—see Section 2.1), we deduce that

$$\begin{aligned} E[Ff_z(F)] &= E[-\delta DL^{-1}Ff_z(F)] \\ &= E[\langle Df_z(F), -DL^{-1}F \rangle_{\mathfrak{H}}] \\ &= E[f'_z(F)\langle DF, -DL^{-1}F \rangle_{\mathfrak{H}}]. \end{aligned}$$

It follows that $E[f'_z(F) - Ff_z(F)] = E[f'_z(F)(1 - \langle DF, -DL^{-1}F \rangle_{\mathfrak{H}})]$ and the proof of the lemma is complete. □

As an application, we deduce the following result, first proven in [20] (the proof is reproduced here for the sake of completeness).

THEOREM 2.4. *Let $F \in \mathbb{D}^{1,2}$ have zero mean and $N \sim \mathcal{N}(0,1)$. Then,*

(2.12) $$d_{\mathrm{Kol}}(F, N) \leq \sqrt{E[(1 - \langle DF, -DL^{-1}F \rangle_{\mathfrak{H}})^2]}.$$

*If $F = I_q(f)$ for some $q \geq 2$ and $f \in \mathfrak{H}^{\odot q}$, then $\langle DF, -DL^{-1}F \rangle_{\mathfrak{H}} = q^{-1} \|DF\|_{\mathfrak{H}}^2$ and therefore*

(2.13) $$d_{\mathrm{Kol}}(F, N) \leq \sqrt{E[(1 - q^{-1}\|DF\|_{\mathfrak{H}}^2)^2]}.$$

PROOF. If $f$ is a bounded, continuously differentiable function such that $\|f'\|_\infty \leq 1$, then, using the same arguments as in the proof of Lemma 2.3 (here, since $f$ belongs to $\mathscr{C}_b^1$, observe that we do not need to assume that the law of $F$ is absolutely continuous), we have

$$\begin{aligned} |E[f'(F) - Ff(F)]| &= |E[f'(F)(1 - \langle DF, -DL^{-1}F \rangle_{\mathfrak{H}})]| \\ &\leq E|1 - \langle DF, -DL^{-1}F \rangle_{\mathfrak{H}}|. \end{aligned}$$

In fact, the inequality $|E[f'(F) - Ff(F)]| \leq E|1 - \langle DF, -DL^{-1}F \rangle_{\mathfrak{H}}|$ continues to hold with $f = f_z$ (which is bounded and Lipschitz, with Lipschitz constant less than one) as is easily seen by convoluting $f_z$ by an approximation of the identity. Hence, Lemma 2.3, combined with Cauchy–Schwarz inequality, implies the desired conclusion. □

REMARK 2.5. In general, the bound appearing on the right-hand side of (2.12) may be infinite. Indeed, the fact that $F \in \mathbb{D}^{1,2}$ only implies that $\langle DF, -DL^{-1}F \rangle_{\mathfrak{H}} \in L^1(\Omega)$. By using the Cauchy–Schwarz inequality twice, one sees that a sufficient condition, in order to have $\langle DF, -DL^{-1}F \rangle_{\mathfrak{H}} \in L^2(\Omega)$, is that $\|DF\|_{\mathfrak{H}}$ and $\|DL^{-1}F\|_{\mathfrak{H}}$ belong to $L^4(\Omega)$. Also, note that if $F$ is equal to a finite sum of multiple integrals (e.g., $F$ is a polynomial functional of $X$), then the random variable $\langle DF, -DL^{-1}F \rangle_{\mathfrak{H}}$ is also a finite sum of multiple integrals and therefore has finite moments of all orders. In particular, for $F = I_q(f)$, the right-hand side of (2.13) is always finite.

The bounds appearing in Theorem 2.4 should be compared with the forthcoming Theorem 2.6, dealing with CLTs on a single Wiener chaos (part A) and on a fixed sum of Wiener chaoses (part B).

THEOREM 2.6 (See [22, 23, 25, 26]). *Fix $q \geq 2$ and let the sequence $F_n = I_q(f_n)$, $n \geq 1$, where $\{f_n\} \subset \mathfrak{H}^{\odot q}$, be such that $E[F_n^2] \to 1$ as $n \to \infty$.*

(A) *The following four conditions are equivalent as $n \to \infty$:*

  (i) $F_n \xrightarrow{\text{Law}} N \sim \mathcal{N}(0,1)$;
  (ii) $E(F_n^4) \to 3$;
  (iii) $\|f_n \otimes_j f_n\|_{\mathfrak{H}^{\otimes 2(q-j)}} \to 0$ *for every* $j = 1, \ldots, q-1$;
  (iv) $1 - q^{-1}\|DF_n\|_{\mathfrak{H}}^2 \to 0$ *in* $L^2$.

(B) *Assume that any one of conditions* (i)–(iv) *of part* A *is satisfied. Let the sequence $G_n$, $n \geq 1$, have the form*

$$G_n = \sum_{p=1}^{M} I_p(g_n^{(p)}), \qquad n \geq 1,$$

*for some $M \geq 1$ (independent of $n$) and some kernels $g_n^{(p)} \in \mathfrak{H}^{\odot p}$ ($p = 1, \ldots, M$, $n \geq 1$). Suppose that, as $n \to \infty$,*

$$E(G_n^2) = \sum_{p=1}^{M} p! \|g_n^{(p)}\|_{\mathfrak{H}^{\otimes p}}^2 \longrightarrow c^2 > 0 \quad \text{and} \quad \|g_n^{(p)} \otimes_j g_n^{(p)}\|_{\mathfrak{H}^{\otimes 2(p-j)}} \longrightarrow 0$$

*for every $p = 1, \ldots, M$ and every $j = 1, \ldots, p-1$. If the sequence of covariances $E(F_n G_n)$ converges to a finite limit, say $\rho \in \mathbb{R}$, then $(F_n, G_n)$ converges in distribution to a two-dimensional Gaussian vector $(N_1, N_2)$ such that $E(N_1^2) = 1$, $E(N_2^2) = c^2$ and $E(N_1 N_2) = \rho$.*

The equivalence between points (i)–(iii) in part A of the previous statement was first proven in [23] by means of stochastic calculus techniques; the fact that condition (iv) is also necessary and sufficient was proven in [22].

Part B (whose proof is straightforward and therefore omitted) is a consequence of the main results established in [25, 26]. Note that in part B of the previous statement, we may allow some of the kernels $g_n^{(p)}$ to be equal to zero. See [19] and [20], Section 3.3, for some extensions of Theorems 2.4 and 2.6 to the framework of noncentral limit theorems.

REMARK ON NOTATION. In what follows, given two numerical sequences $\{a_n\}$ and $\{b_n\}$, the symbol $a_n \sim b_n$ means that $\lim a_n/b_n = 1$, whereas $a_n \asymp b_n$ means that the ratio $a_n/b_n$ converges to a nonzero finite limit.

2.3. *A useful computation.* We shall denote by $\{H_q : q \geq 0\}$ the class of Hermite polynomials, defined as follows: $H_0 \equiv 1$ and, for $q \geq 1$,

$$(2.14) \qquad H_q(z) = (-1)^q e^{z^2/2} \frac{d^q}{dz^q} e^{-z^2/2}, \qquad z \in \mathbb{R};$$

for instance, $H_1(z) = z$, $H_2(z) = z^2 - 1$ and so on. Note that the definition of the class $\{H_q\}$ immediately implies the recurrence relation

$$(2.15) \qquad \frac{d}{dz} H_q(z) e^{-z^2/2} = -H_{q+1}(z) e^{-z^2/2},$$

yielding that the Hermite polynomials are related to the derivatives of $\Phi(z) = P(N \leq z)$ [$N \sim \mathcal{N}(0,1)$], written $\Phi^{(q)}(z)$ ($q = 1, 2, \ldots$), by the formula

$$(2.16) \qquad \Phi^{(q)}(z) = (-1)^{q-1} H_{q-1}(z) \frac{e^{-z^2/2}}{\sqrt{2\pi}}.$$

We also have, for any $q \geq 1$,

$$(2.17) \qquad \frac{d}{dz} H_q(z) = q H_{q-1}(z).$$

Now, denote by $f_z$ the solution to the Stein equation (1.5) given in formulae (2.10)–(2.11). The following result, connecting $f_z$ with the Hermite polynomials and the derivatives of $\Phi$, will be used in Section 3.

PROPOSITION 2.7. *For every $q \geq 1$ and every $z \in \mathbb{R}$,*

$$(2.18) \qquad \begin{aligned} \int_{-\infty}^{+\infty} f_z'(x) H_q(x) \frac{e^{-x^2/2}}{\sqrt{2\pi}} \, dx &= \frac{1}{q+2} H_{q+1}(z) \frac{e^{-z^2/2}}{\sqrt{2\pi}} \\ &= \frac{1}{q+2} (-1)^{q+1} \Phi^{(q+2)}(z). \end{aligned}$$

PROOF. By integrating by parts and by exploiting relations (2.11) and (2.15), one obtains that

$$\int_{-\infty}^{+\infty} f'_z(x) H_q(x) \frac{e^{-x^2/2}}{\sqrt{2\pi}} \, dx$$

(2.19)
$$= \int_{-\infty}^{+\infty} f_z(x) H_{q+1}(x) \frac{e^{-x^2/2}}{\sqrt{2\pi}} \, dx$$

$$= \frac{1}{\sqrt{2\pi}} \int_{-\infty}^{+\infty} H_{q+1}(x) \left( \int_{-\infty}^{x} (\mathbf{1}_{(-\infty,z]}(a) - \Phi(z)) e^{-a^2/2} \, da \right) dx.$$

By integrating by parts, using $H_{q+1} = \frac{1}{q+2} H'_{q+2}$ [see (2.17)] and in view of (2.15), one easily proves that

$$\int_{-\infty}^{+\infty} H_{q+1}(x) \left( \int_{-\infty}^{x} (\mathbf{1}_{(-\infty,z]}(a) - \Phi(z)) e^{-a^2/2} \, da \right) dx$$

$$= -\frac{1}{q+2} \int_{-\infty}^{+\infty} H_{q+2}(x)(\mathbf{1}_{(-\infty,z]}(x) - \Phi(z)) e^{-x^2/2} \, dx$$

$$= -\frac{1}{q+2} \left( \int_{-\infty}^{z} H_{q+2}(x) e^{-x^2/2} \, dx - \Phi(z) \int_{-\infty}^{+\infty} H_{q+2}(x) e^{-x^2/2} \, dx \right)$$

$$= \frac{1}{q+2} H_{q+1}(z) e^{-z^2/2}.$$

By plugging this expression into (2.19), we immediately arrive at the desired conclusion. □

For instance, by specializing formula (2.18) to the case $q = 1$, one obtains, for $N \sim \mathcal{N}(0,1)$,

(2.20) $$E[f'_z(N) \times N] = \frac{1}{3}(z^2 - 1) \frac{e^{-z^2/2}}{\sqrt{2\pi}} = \frac{1}{3} \Phi^{(3)}(z).$$

## 3. Main results.

3.1. *Two general statements.* We start by studying the case of a general sequence of Malliavin derivable functionals.

THEOREM 3.1. *Let $F_n$, $n \geq 1$, be a sequence of centered and square-integrable functionals of some isonormal Gaussian process $X = \{X(h) : h \in \mathfrak{H}\}$ such that $E(F_n^2) \longrightarrow 1$ as $n \to \infty$. Suppose that the following three conditions hold:*

(i) *for every $n$, one has that $F_n \in \mathbb{D}^{1,2}$ and $F_n$ has an absolutely continuous law (with respect to the Lebesgue measure);*

(ii) *the quantity* $\varphi(n) = \sqrt{E[(1 - \langle DF_n, -DL^{-1}F_n\rangle_{\mathfrak{H}})^2]}$ *[as in (1.8)] is such that:* (a) $\varphi(n)$ *is finite for every* $n$; (b) *as* $n \to \infty$, $\varphi(n)$ *converges to zero; and* (c) *there exists* $m \geq 1$ *such that* $\varphi(n) > 0$ *for* $n \geq m$;

(iii) *as* $n \to \infty$, *the two-dimensional vector* $(F_n, \frac{1 - \langle DF_n, -DL^{-1}F_n\rangle_{\mathfrak{H}}}{\varphi(n)})$ *[as in formula (1.10)] converges in distribution to a centered two-dimensional Gaussian vector* $(N_1, N_2)$ *such that* $E(N_1^2) = E(N_2^2) = 1$ *and* $E(N_1 N_2) = \rho$.

*Then, the upper bound* $d_{\mathrm{Kol}}(F_n, N) \leq \varphi(n)$ *holds. Moreover, for every* $z \in \mathbb{R}$,

$$(3.1) \quad \varphi(n)^{-1}[P(F_n \leq z) - \Phi(z)] \xrightarrow[n \to \infty]{} \frac{\rho}{3}(z^2 - 1)\frac{e^{-z^2/2}}{\sqrt{2\pi}} = \frac{\rho}{3}\Phi^{(3)}(z).$$

*As a consequence, if* $\rho \neq 0$, *then there exists a constant* $c \in (0,1)$, *as well as an integer* $n_0 \geq 1$, *such that relation (1.3) holds for* $d = d_{\mathrm{Kol}}$ *and for every* $n \geq n_0$.

PROOF.  Fix $z \in \mathbb{R}$. From assumption (i) and Lemma 2.3, recall that

$$\varphi(n)^{-1}[P(F_n \leq z) - \Phi(z)] = E[f'_z(F_n)\varphi(n)^{-1}(1 - \langle DF_n, -DL^{-1}F_n\rangle_{\mathfrak{H}})].$$

The facts that $f'_z$ is bounded by 1 on the one hand and that $\varphi(n)^{-1}(1 - \langle DF_n, -DL^{-1}F_n\rangle_{\mathfrak{H}})$ has variance 1 on the other hand imply that the sequence

$$f'_z(F_n)\varphi(n)^{-1}(1 - \langle DF_n, -DL^{-1}F_n\rangle_{\mathfrak{H}}), \qquad n \geq 1,$$

is uniformly integrable. Now, deduce from (2.10) that $x \to f'_z(x)$ is continuous at every $x \neq z$. This yields that, as $n \to \infty$ and due to assumption (iii),

$$E[f'_z(F_n)\varphi(n)^{-1}(1 - \langle DF_n, -DL^{-1}F_n\rangle_{\mathfrak{H}})]$$
$$\longrightarrow E(f'_z(N_1)N_2) = \rho \times E(f'_z(N_1)N_1).$$

Relation (3.1) now follows from formula (2.20). If, in addition, $\rho \neq 0$, then one can obtain the lower bound (1.3) by using the elementary relation

$$|P(F_n \leq 0) - \Phi(0)| \leq d_{\mathrm{Kol}}(F_n, N). \qquad \square$$

REMARK 3.2.  Clearly, the conclusion of Theorem 3.1 still holds when $n$ is replaced by some continuous parameter. The same remark holds for the forthcoming results of this section.

The next proposition connects our results with one-term Edgeworth expansions. Note that, in the following statement, we assume that $E(F_n) = 0$ and $E(F_n^2) = 1$ so that the first term in the (formal) Edgeworth expansion of $P(F_n \leq z) - \Phi(z)$ coincides with $-(3!)^{-1}E(F_n^3)\Phi^{(3)}(z)$. For an introduction to Edgeworth expansions, the reader is referred, for example, to McCullagh [18], Chapter 3, or Hall [12], Chapter 2. See also Rotar [31] for another application of Stein's method to Edgeworth expansions.

PROPOSITION 3.3 (One-term Edgeworth expansions). *Let $F_n$, $n \geq 1$, be a sequence of centered and square-integrable functionals of the isonormal Gaussian process $X = \{X(h) : h \in \mathfrak{H}\}$ such that $E(F_n^2) = 1$. Suppose that conditions* (i)–(iii) *of Theorem 3.1 are satisfied and also that:*

(a) *for every $n$, $E|F_n|^3 < \infty$;*
(b) *there exists $\varepsilon > 0$ such that $\sup_{n \geq 1} E|F_n|^{2+\varepsilon} < \infty$.*

*Then, as $n \to \infty$,*

$$(3.2) \qquad \frac{1}{2\varphi(n)} E(F_n^3) \longrightarrow -\rho$$

*and, for every $z \in \mathbb{R}$, one has the one-term local Edgeworth expansion*

$$(3.3) \qquad P(F_n \leq z) - \Phi(z) + \frac{1}{3!} E(F_n^3) \Phi^{(3)}(z) = o_z(\varphi(n)),$$

*where $o_z(\varphi(n))$ indicates a numerical sequence (depending on $z$) such that $\varphi(n)^{-1} o_z(\varphi(n)) \to 0$ as $n \to \infty$.*

REMARK 3.4. Of course, relation (3.3) is interesting only when $\rho \neq 0$. Indeed, in this case, one has that, thanks to Theorem 3.1, $P(F_n \leq z) - \Phi(z) \asymp \varphi(n)$ (the symbol $\asymp$ indicates asymptotic equivalence) so that, for a fixed $z$, the addition of $\frac{1}{3!} E(F_n^3) \Phi^{(3)}(z)$ actually increases the rate of convergence to zero.

PROOF OF PROPOSITION 3.3. Since assumption (a) is in order and $E(F_n) = 0$, one can deduce from Lemma 2.1, in the case $s = 1$, that

$$E\left(F_n \times \frac{1 - \langle DF_n, -DL^{-1}F_n \rangle_{\mathfrak{H}}}{\varphi(n)}\right) = -\frac{1}{2\varphi(n)} E(F_n^3).$$

Assumption (b), combined with the fact that $\varphi(n)^{-1}(1 - \langle DF_n, -DL^{-1}F_n \rangle_{\mathfrak{H}})$ has variance 1, immediately yields that there exists $\delta > 0$ such that

$$\sup_{n \geq 1} E|F_n \times \varphi(n)^{-1}(1 - \langle DF_n, -DL^{-1}F_n \rangle_{\mathfrak{H}})|^{1+\delta} < \infty.$$

In particular, the sequence $\{F_n \times \varphi(n)^{-1}(1 - \langle DF_n, -DL^{-1}F_n \rangle_{\mathfrak{H}}) : n \geq 1\}$ is uniformly integrable. Therefore, since assumption (iii) in the statement of Theorem 3.1 is in order, we may deduce that, as $n \to \infty$,

$$\frac{1}{2\varphi(n)} E(F_n^3) \longrightarrow -E(N_1 N_2) = -\rho.$$

As a consequence,

$$\varphi(n)^{-1} \left| P(F_n \leq z) - \Phi(z) + \frac{1}{3!} E(F_n^3) \Phi^{(3)}(z) \right|$$

$$\leq \left| \frac{P(F_n \leq z) - \Phi(z)}{\varphi(n)} - \frac{\rho}{3} \Phi^{(3)}(z) \right| + \frac{|\Phi^{(3)}(z)|}{3} \left| \frac{1}{2\varphi(n)} E(F_n^3) + \rho \right|$$

and the conclusion follows from Theorem 3.1. □

REMARK 3.5. By inspection of the proof of Proposition 3.3, one sees that Assumption (b) in the statement may equally well be replaced by the following, weaker, condition: (b′) *the sequence*

$$F_n \times \varphi(n)^{-1}(1 - \langle DF_n, -DL^{-1}F_n \rangle_{\mathfrak{H}}), \qquad n \geq 1,$$

*is uniformly integrable.*

3.2. *Multiple integrals.* The following statement specializes the content of the previous subsection to multiple integrals with respect to some isonormal Gaussian process $X = \{X(h) : h \in \mathfrak{H}\}$. Recall that a nonzero finite sum of multiple integrals of arbitrary orders is always an element of $\mathbb{D}^{1,2}$ and, also, that its law admits a density with respect to Lebesgue measure [this implies that assumption (i) in the statement of Theorem 3.1 is automatically satisfied in this section]; see Shigekawa [32].

PROPOSITION 3.6. *Let $q \geq 2$ be a fixed integer and let the sequence $F_n$, $n \geq 1$, have the form $F_n = I_q(f_n)$, where, for $n \geq 1$, $f_n \in \mathfrak{H}^{\odot q}$. Suppose that $E(F_n^2) = q!\|f_n\|_{\mathfrak{H}^{\otimes q}}^2 \to 1$. Then, the quantity $\varphi(n)$ appearing in formula (1.8) is such that*

$$(3.4) \quad \varphi(n)^2 = E[(1 - q^{-1}\|DF_n\|_{\mathfrak{H}}^2)^2]$$

$$(3.5) \qquad = (1 - q!\|f_n\|_{\mathfrak{H}^{\otimes q}}^2)^2$$

$$+ q^2 \sum_{r=1}^{q-1} (2q - 2r)!(r-1)!^2 \binom{q-1}{r-1}^4 \|f_n \widetilde{\otimes}_r f_n\|_{\mathfrak{H}^{\otimes 2(q-r)}}^2.$$

*Now, suppose that, as $n \to \infty$,*

$$(3.6) \qquad \|f_n \otimes_r f_n\|_{\mathfrak{H}^{\otimes 2(q-r)}} \to 0$$

*for every $r = 1, \ldots, q-1$ and, also,*

$$(3.7) \qquad \frac{1 - q!\|f_n\|_{\mathfrak{H}^{\otimes q}}^2}{\varphi(n)} \longrightarrow 0.$$

*Then, assumption* (ii) *in the statement of Theorem 3.1 is satisfied and a set of sufficient conditions, implying that assumption* (iii) *in the same theorem holds, are the following relations (3.8)–(3.9): as $n \to \infty$,*

$$(3.8) \qquad \varphi(n)^{-2}\|(f_n \widetilde{\otimes}_r f_n) \otimes_l (f_n \widetilde{\otimes}_r f_n)\|_{\mathfrak{H}^{\otimes 2(2(q-r)-l)}} \to 0$$

*for every $r = 1, \ldots, q-1$ and every $l = 1, \ldots, 2(q-r) - 1$ and, if $q$ is even,*

$$(3.9) \quad -qq!(q/2-1)!\binom{q-1}{q/2-1}^2 \varphi(n)^{-1} \langle f_n, f_n \widetilde{\otimes}_{q/2} f_n \rangle_{\mathfrak{H}^{\otimes q}} \longrightarrow \rho.$$

*If $q$ is odd and (3.8) holds, then assumption* (ii) *in Theorem 3.1 holds with $\rho = 0$.*

PROOF. Formulae (3.4)–(3.5) are a consequence of [20], Proposition 3.2. The fact that (3.6) implies $\varphi(n) \longrightarrow 0$ is immediate (recall that $\|f_n \otimes_r f_n\|_{\mathfrak{H}^{\otimes 2(q-r)}} \geq \|f_n \widetilde{\otimes}_r f_n\|_{\mathfrak{H}^{\otimes 2(q-r)}}$). Again using [20], formula (3.42), one has that

$$
\begin{aligned}
(3.10) \quad & \frac{1 - q^{-1}\|DI_q(f_n)\|_{\mathfrak{H}}^2}{\varphi(n)} \\
& = \frac{1 - q!\|f_n\|_{\mathfrak{H}^{\otimes q}}^2}{\varphi(n)} - q \sum_{r=1}^{q-1} (r-1)! \binom{q-1}{r-1}^2 I_{2(q-r)} \left( \frac{f_n \widetilde{\otimes}_r f_n}{\varphi(n)} \right).
\end{aligned}
$$

Finally, the fact that (3.8) and (3.9) (for $q$ even) imply that assumption (iii) in Theorem 3.1 is satisfied, is a consequence of representation (3.10) and part B of Theorem 2.6, in the case

$$
G_n = -q \sum_{r=1}^{q-1} (r-1)! \binom{q-1}{r-1}^2 I_{2(q-r)} \left( \frac{f_n \widetilde{\otimes}_r f_n}{\varphi(n)} \right),
$$

and $c^2 = 1$, by taking into account the fact that, for $q$ even,

$$
E(F_n G_n) = -q q! (q/2 - 1)! \binom{q-1}{q/2-1}^2 \varphi(n)^{-1} \langle f_n, f_n \widetilde{\otimes}_{q/2} f_n \rangle_{\mathfrak{H}^{\otimes q}},
$$

whereas $E(F_n G_n) = 0$ for $q$ odd. $\square$

REMARK 3.7. Observe that, due to part A of Theorem 2.6, condition (3.6) is actually necessary and sufficient to have $\varphi(n) \longrightarrow 0$. Moreover, if conditions (3.6)–(3.9) are satisfied, then the usual properties of finite sums of multiple integrals (see, e.g., [13], Chapter VI) imply that assumptions (a)–(b) in the statement of Proposition 3.3 are automatically met so that Proposition 3.6 indeed provides one-term local Edgeworth expansions.

3.3. *Second Wiener chaos.* In this subsection, we focus on random variables in the second Wiener chaos associated with an isonormal Gaussian process $X = \{X(h) : h \in \mathfrak{H}\}$, that is, random variables of the type $F = I_2(f)$, where $f \in \mathfrak{H}^{\odot 2}$. To every kernel $f \in \mathfrak{H}^{\odot 2}$, we associate two objects: (I) the Hilbert–Schmidt operator

$$
(3.11) \qquad H_f : \mathfrak{H} \to \mathfrak{H}; \qquad g \to f \otimes_1 g,
$$

where the contraction $f \otimes_1 g$ is defined according to (2.3), and (II) the sequence of kernels $\{f \otimes_1^{(p)} f : p \geq 1\} \subset \mathfrak{H}^{\odot 2}$, defined as follows: $f \otimes_1^{(1)} f = f$ and, for $p \geq 2$,

$$
(3.12) \qquad f \otimes_1^{(p)} f = (f \otimes_1^{(p-1)} f) \otimes_1 f.
$$

We write $\{\lambda_{f,j}\}_{j\geq 1}$ to indicate the eigenvalues of $H_f$. Now, for $p \geq 1$, denote by $\kappa_p(I_2(f))$ the $p$th cumulant of $I_2(f)$. The following relation, giving an explicit expression for the cumulants of $I_2(f)$, is well known (see, e.g., [6] for a proof): one has that $\kappa_1(I_2(f)) = E(I_2(f)) = 0$ and, for $p \geq 2$,

$$
\begin{aligned}
\kappa_p(I_2(f)) &= 2^{p-1}(p-1)! \times \mathrm{Tr}(H_f^p) \\
&= 2^{p-1}(p-1)! \times \langle f \otimes_1^{(p-1)} f, f \rangle_{\mathfrak{H}^{\otimes 2}} \\
&= 2^{p-1}(p-1)! \times \sum_{j=1}^{\infty} \lambda_{f,j}^p,
\end{aligned}
\tag{3.13}
$$

where $\mathrm{Tr}(H_f^p)$ stands for the trace of the $p$th power of $H_f$.

PROPOSITION 3.8. *Let $F_n = I_2(f_n)$, $n \geq 1$, be such that $f_n \in \mathfrak{H}^{\odot 2}$ and write $\kappa_p^{(n)} = \kappa_p(F_n)$, $p \geq 1$. Assume that $\kappa_2^{(n)} = E(F_n^2) \longrightarrow 1$ as $n \to \infty$. Then, as $n \to \infty$, $F_n \xrightarrow{\mathrm{Law}} N \sim \mathcal{N}(0,1)$ if and only if $\kappa_4^{(n)} \longrightarrow 0$. In this case, we further have*

$$
d_{\mathrm{Kol}}(F_n, N) \leq \sqrt{\frac{\kappa_4^{(n)}}{6} + (\kappa_2^{(n)} - 1)^2}.
\tag{3.14}
$$

*If, in addition, we have, as $n \to \infty$,*

$$
\frac{\kappa_2^{(n)} - 1}{\kappa_4^{(n)}/6 + (\kappa_2^{(n)} - 1)^2} \longrightarrow 0,
\tag{3.15}
$$

$$
\frac{\kappa_3^{(n)}}{\sqrt{\kappa_4^{(n)}/6 + (\kappa_2^{(n)} - 1)^2}} \longrightarrow \alpha \quad and \quad \frac{\kappa_8^{(n)}}{(\kappa_4^{(n)}/6 + (\kappa_2^{(n)} - 1)^2)^2} \longrightarrow 0,
\tag{3.16}
$$

*then*

$$
\frac{P(F_n \leq z) - \Phi(z)}{\sqrt{\kappa_4^{(n)}/6 + (\kappa_2^{(n)} - 1)^2}} \longrightarrow \frac{\alpha}{3!} \frac{1}{\sqrt{2\pi}} (1 - z^2) e^{-z^2/2} \qquad as\ n \to \infty.
\tag{3.17}
$$

*In particular, if $\alpha \neq 0$, then there exists $c \in (0,1)$ and $n_0 \geq 1$ such that, for any $n \geq n_0$,*

$$
\sup_{z \in \mathbb{R}} |P(F_n \leq z) - \Phi(z)| \geq c \sqrt{\frac{\kappa_4^{(n)}}{6} + (\kappa_2^{(n)} - 1)^2}.
\tag{3.18}
$$

REMARK 3.9. 1. If $E(F_n^2) = \kappa_2^{(n)} = 1$, then condition (3.15) becomes immaterial and the denominators appearing in formula (3.16) involve solely $\kappa_4^{(n)}$.

2. By combining (3.16) with (3.17), we have that, as $n \to \infty$,

$$P(F_n \leq z) - \Phi(z) \sim \frac{\kappa_3^{(n)}}{3!\sqrt{2\pi}}(1-z^2)e^{-z^2/2},$$

whenever $z \neq \pm 1$ and $\alpha \neq 0$.

PROOF OF PROPOSITION 3.8. First, since $E(F_n) = 0$, we have $\kappa_4^{(n)} = E(F_n^4) - 3E(F_n^2)^2$. Thus, the equivalence between $\kappa_4^{(n)} \longrightarrow 0$ and $F_n \xrightarrow{\text{Law}} \mathcal{N}(0,1)$ is a direct consequence of part A of Theorem 2.6. Now, observe that

$$\tfrac{1}{2}\|DF_n\|^2 - 1 = 2I_2(f_n \otimes_1 f_n) + E(F_n^2) - 1 = 2I_2(f_n \otimes_1 f_n) + \kappa_2^{(n)} - 1.$$

In particular,

$$\text{Var}\left(\frac{1}{2}\|DZ_n\|^2 - 1\right) = 8\|f_n \otimes_1 f_n\|_{\mathfrak{H}^{\otimes 2}}^2 + (\kappa_2^{(n)} - 1)^2 = \frac{\kappa_4^{(n)}}{6} + (\kappa_2^{(n)} - 1)^2,$$

where we have used (3.13) in the case $p = 4$ (note that $\langle f \otimes_1^{(3)} f, f\rangle_{\mathfrak{H}^{\otimes 2}} = \|f \otimes_1 f\|_{\mathfrak{H}^{\otimes 2}}^2$). This implies that the quantity $\varphi(n)$ appearing in (1.7) indeed equals $\sqrt{\kappa_4^{(n)}/6 + (\kappa_2^{(n)} - 1)^2}$. To conclude the proof, it is sufficient to apply Proposition 3.6 in the case $q = 2$, by observing that

$$\frac{1 - \kappa_2^{(n)}}{\sqrt{\kappa_4^{(n)}/6 + (\kappa_2^{(n)} - 1)^2}} = \frac{1 - 2\|f_n\|_{\mathfrak{H}^{\otimes 2}}^2}{\varphi(n)}$$

and, also, by using (3.13) in the cases $p = 3$ and $p = 8$, respectively,

$$\frac{\kappa_3^{(n)}}{\sqrt{\kappa_4^{(n)}/6 + (\kappa_2^{(n)} - 1)^2}} = \frac{8\langle f_n, f_n \otimes_1 f_n\rangle_{\mathfrak{H}^{\otimes 2}}}{\varphi(n)}$$

and

$$\frac{\kappa_8^{(n)}}{(\kappa_4^{(n)}/6 + (\kappa_2^{(n)} - 1)^2)^2} = 2^7 7! \times \frac{\langle f \otimes_1^{(7)} f, f\rangle_{\mathfrak{H}^{\otimes 2}}}{\varphi(n)^4}$$

$$= 2^7 7! \times \frac{\|(f_n \otimes_1 f_n) \otimes_1 (f_n \otimes_1 f_n)\|_{\mathfrak{H}^{\otimes 2}}^2}{\varphi(n)^4}. \qquad \square$$

**4. Toeplitz quadratic functionals of continuous-time stationary processes.** In this section, we apply our results to establish (possibly optimal) Berry–Esseen bounds in CLTs involving quadratic functionals of continuous-time stationary Gaussian processes. Our results represent a substantial refinement of the CLTs proven in the papers by Ginovyan [7] and Ginovyan and Sahakyan [8], where the authors have extended to a continuous-time setting the discrete-time results of Avram [1], Fox and Taqqu [6] and Giraitis and Surgailis [10]. In the discrete-time case, Berry–Esseen-type bounds for CLTs involving special quadratic functionals of stationary Gaussian processes are obtained in [35], and Edgeworth expansions are studied in, for example, [17]. However, to our knowledge, the results which are proved in this section are the first (exact) Berry–Esseen bounds ever proved in the continuous-time case. Observe that it is not clear whether one can deduce bounds in continuous-time by using the discrete-time findings of [17] and [35]. We refer the reader to [2] and [11] (and the references therein) for CLTs and one-term Edgeworth expansions concerning quadratic functionals of general discrete-time processes.

Let $X = (X_t)_{t \in \mathbb{R}}$ be a centered real-valued Gaussian process with spectral density $f : \mathbb{R} \to \mathbb{R}$. This means that, for every $u, t \in \mathbb{R}$, one has

$$E(X_u X_{u+t}) := r(t) = \int_{-\infty}^{+\infty} e^{i\lambda t} f(\lambda) \, d\lambda, \qquad t \in \mathbb{R},$$

where $r : \mathbb{R} \to \mathbb{R}$ is the covariance function of $X$. We stress that the density $f$ is necessarily an even function. For $T > 0$, let $Q_T = \iint_{[0,T]^2} \widehat{g}(t-s) X_t X_s \, dt \, ds$, where

$$\widehat{g}(t) = \int_{-\infty}^{+\infty} e^{i\lambda t} g(\lambda) \, d\lambda, \qquad t \in \mathbb{R},$$

is the Fourier transform of some integrable *even* function $g : \mathbb{R} \to \mathbb{R}$. The random variable $Q_T$ is customarily called the *Toeplitz quadratic functional* of $X$, associated with $g$ and $T$. We also set

$$\tilde{Q}_T = \frac{Q_T - E(Q_T)}{\sqrt{T}} \quad \text{and} \quad \check{Q}_T = \frac{\tilde{Q}_T}{\sigma(T)}$$

with $\sigma(T)^2 = \text{Var}(\tilde{Q}_T)$. The cumulants of $\tilde{Q}_T$ and $\check{Q}_T$ are denoted, respectively, by

$$\tilde{\kappa}_j^{(T)} = \kappa_j(\tilde{Q}_T) \quad \text{and} \quad \check{\kappa}_j^{(T)} = \check{\kappa}_j(\tilde{Q}_T), \qquad j \geq 1.$$

Given $T > 0$ and $\psi \in L^1(\mathbb{R})$, we denote by $B_T(\psi)$ the *truncated Toeplitz operator* associated with $\psi$ and $T$, acting on a square-integrable function $u$ as follows:

$$B_T(\psi)(u)(\lambda) = \int_0^T u(x) \widehat{\psi}(\lambda - x) \, dx, \qquad \lambda \in \mathbb{R},$$

where $\widehat{\psi}$ is the Fourier transform of $\psi$. Given $\psi, \gamma \in L^1(\mathbb{R})$, we denote by $B_T(\psi) B_T(\gamma)$ the product of the two operators $B_T(\psi)$ and $B_T(\gamma)$; also, $[B_T(\psi) B_T(\gamma)]^j$, $j \geq 1$, is the $j$th power of $B_T(\psi) B_T(\gamma)$. The symbol $\text{Tr}(U)$ indicates the trace of an operator $U$.

The following statement collects some of the results proven in [7, 8].

THEOREM 4.1 (See [7, 8]).  1. *For every $j \geq 1$, the $j$th cumulant of $\tilde{Q}_T$ is given by*

$$\tilde{\kappa}_j^{(T)} = \begin{cases} 0, & \text{if } j = 1, \\ T^{-j/2} 2^{j-1}(j-1)! \operatorname{Tr}[B_T(f)B_T(g)]^j, & \text{if } j \geq 2. \end{cases}$$

2. *Assume that $f \in L^p(\mathbb{R}) \cap L^1(\mathbb{R})$ ($p \geq 1$), that $g \in L^q(\mathbb{R}) \cap L^1(\mathbb{R})$ ($q \geq 1$) and that $\frac{1}{p} + \frac{1}{q} \leq \frac{1}{j}$. Then,*

$$\tilde{\kappa}_j^{(T)} \underset{T \to \infty}{\sim} T^{1-j/2} \times 2^{j-1}(j-1)!(2\pi)^{2j-1} \int_{-\infty}^{+\infty} f(x)^j g(x)^j \, dx.$$

3. *If $\frac{1}{p} + \frac{1}{q} \leq \frac{1}{2}$, then*

$$\sigma^2(T) = \tilde{\kappa}_2^{(T)} \underset{T \to \infty}{\longrightarrow} \sigma^2(\infty) := 16\pi^3 \int_{-\infty}^{+\infty} f^2(x) g^2(x) \, dx$$

*and $\check{Q}_T \overset{\text{Law}}{\longrightarrow} Z \sim \mathcal{N}(0, 1)$ as $T \to \infty$.*

The next statement shows that one can apply Proposition 3.8 in order to obtain Berry–Esseen bounds for the CLT appearing in point 3 of Theorem 4.1. Observe that, since the variance of $\check{Q}_T$ is equal to 1, by construction, to establish an upper bound, we need to control only the fourth cumulant of $\check{Q}_T$: this will be done by using point 2 of Theorem 4.1 and by assuming that $\frac{1}{p} + \frac{1}{q} \leq \frac{1}{4}$. On the other hand, to prove lower bounds, one needs to have a precise estimate of the asymptotic behavior of the eighth cumulant of $\check{Q}_T$: again in view of point 2 of Theorem 4.1, this requires that $\frac{1}{p} + \frac{1}{q} \leq \frac{1}{8}$.

THEOREM 4.2. *Assume that $f \in L^p(\mathbb{R}) \cap L^1(\mathbb{R})$ ($p \geq 1$) and that $g \in L^q(\mathbb{R}) \cap L^1(\mathbb{R})$ ($q \geq 1$). Let $\Phi(z) = P(N \leq z)$, where $N \sim \mathcal{N}(0, 1)$.*

1. *If $\frac{1}{p} + \frac{1}{q} \leq \frac{1}{4}$, then there exists $C = C(f, g) > 0$ such that, for all $T > 0$, we have*

$$\sup_{z \in \mathbb{R}} |P(\check{Q}_T \leq z) - \Phi(z)| \leq \frac{C}{\sqrt{T}}.$$

2. *If $\frac{1}{p} + \frac{1}{q} \leq \frac{1}{8}$ and*

$$\int_{-\infty}^{+\infty} f^3(x) g^3(x) \, dx \neq 0,$$

*then there exists $c = c(f, g) > 0$ and $T_0 = T_0(f, g)$ such that $T \geq T_0$ implies*

$$\sup_{z \in \mathbb{R}} |P(\check{Q}_T \leq z) - \Phi(z)| \geq \frac{c}{\sqrt{T}}.$$

*More precisely, for any $z \in \mathbb{R}$, we have*

$$\sqrt{T}(P(\check{Q}_T \leq z) - \Phi(z)) \tag{4.1}$$
$$\xrightarrow[T \to \infty]{} \sqrt{\frac{2}{3}} \frac{\int_{-\infty}^{+\infty} f^3(x) g^3(x)\, dx}{(\int_{-\infty}^{+\infty} f^2(x) g^2(x)\, dx)^{3/2}} (1-z^2) e^{-z^2/2}.$$

PROOF. It is a standard result that each random variable $\check{Q}_T$ can be represented as a double Wiener–Itô integral with respect to $X$. It follows that the statement can be proven by means of Proposition 3.8. Now, whenever $\frac{1}{p} + \frac{1}{q} \leq \frac{1}{j}$, one easily obtains, from points 2 and 3 of Theorem 4.1, that

$$\check{\kappa}_j^{(T)} \underset{T \to \infty}{\sim} T^{1-j/2} \frac{2^{j-1}(j-1)!(2\pi)^{2j-1}}{(16\pi^3)^{j/2}} \frac{\int_{-\infty}^{+\infty} f^j(x) g^j(x)\, dx}{(\int_{-\infty}^{+\infty} f^2(x) g^2(x)\, dx)^{j/2}} \tag{4.2}$$

and the desired conclusion is then obtained by a direct application of Proposition 3.8. In particular, point 1 in the statement is immediately deduced from the fact that $\frac{1}{p} + \frac{1}{q} \leq \frac{1}{4}$, from relation (4.2) and the bound (3.14), with $\check{\kappa}_4^{(T)}$ replacing $\kappa_4^{(n)}$ (observe that $\check{\kappa}_2^{(T)} = 1$, by construction). On the other hand, point 2 is a consequence of the fact that if $\frac{1}{p} + \frac{1}{q} \leq \frac{1}{8}$, then (4.2) implies that condition (3.16) is met. The exact value of the constant appearing on the right-hand side of (4.1) is deduced from elementary simplifications. □

**5. Exploding quadratic functionals of a Brownian sheet.** In this section, we apply our results to the study of some quadratic functionals of a standard Brownian sheet on $[0,1]^d$ ($d \geq 1$), denoted $\mathbf{W} = \{\mathbf{W}(t_1, \ldots, t_d) : (t_1, \ldots, t_d) \in [0,1]^d\}$. We recall that $\mathbf{W}$ is a centered Gaussian process such that, for every $(t_1, \ldots, t_d), (u_1, \ldots, u_d) \in [0,1]^d$,

$$E[\mathbf{W}(t_1, \ldots, t_d) \mathbf{W}(u_1, \ldots, u_d)] = \prod_{i=1,\ldots,d} \min(u_i, t_i)$$

so that, if $d=1$, the process $\mathbf{W}$ is indeed a standard Brownian motion on $[0,1]$. It is easily proved that, for every $d \geq 1$, the Gaussian space generated by $\mathbf{W}$ can be identified with an isonormal Gaussian process of the type $X = \{X(f) : f \in L^2([0,1]^d, d\lambda)\}$, where $d\lambda$ indicates the restriction of Lebesgue measure on $[0,1]^d$. It is also well known that the trajectories of $\mathbf{W}$ enjoy the following, remarkable, property:

$$\int_{[0,1]^d} \left( \frac{\mathbf{W}(t_1, \ldots, t_d)}{t_1 \cdots t_d} \right)^2 dt_1 \cdots dt_d = +\infty, \qquad P\text{-a.s.} \tag{5.1}$$

Relation (5.1) is a consequence of the scaling properties of $\mathbf{W}$ and of the well known *Jeulin's lemma* (see [14], Lemma 1, page 44, or [24]). In the case $d=1$, the study of phenomena such as (5.1) arose at the end of the 1970s, in connection with the theory of enlargement of filtrations (see [14, 15]); see also [16] for some relations with noncanonical representations of Gaussian processes.

Now, denote, for every $\varepsilon > 0$,

$$B_\varepsilon^d = \left\{\int_{[\varepsilon,1]^d} \left(\frac{\mathbf{W}(t_1,\ldots,t_d)}{t_1 \cdots t_d}\right)^2 dt_1 \cdots dt_d\right\} - \left(\log \frac{1}{\varepsilon}\right)^d$$

and observe that $B_\varepsilon^d$ is a centered random variable with moments of all orders. The CLT stated in the forthcoming proposition gives some insights into the "rate of explosion around zero" of the random function

$$(t_1,\ldots,t_d) \to \left(\frac{\mathbf{W}(t_1,\ldots,t_d)}{t_1 \cdots t_d}\right)^2.$$

PROPOSITION 5.1. *For every $d \geq 1$, as $\varepsilon \to 0$,*

(5.2) $$\widetilde{B}_\varepsilon^d := (4\log 1/\varepsilon)^{-d/2} \times B_\varepsilon^d \xrightarrow{\text{Law}} N \sim \mathcal{N}(0,1).$$

Proposition 5.1 has been established in [27] (for the case $d=1$), [5] (for the case $d=2$) and [23] (for the case $d>2$). See [27, 28] for an application of the CLT (5.2) (in the case $d=1$) to the study of Brownian local times. See [5] for some applications to conditioned bivariate Gaussian processes and to statistical tests of independence. The next result, which is obtained by means of the techniques developed in this paper, gives an exact description (in terms of the Kolmogorov distance) of the rate of convergence of $\widetilde{B}_\varepsilon^d$ toward a Gaussian random variable.

PROPOSITION 5.2. *For every $d \geq 1$, there exist constants $0 < c(d) < C(d) < +\infty$ and $0 < \eta(d) < 1$, depending uniquely on $d$, such that, for every $\varepsilon > 0$,*

$$d_{\text{Kol}}[\widetilde{B}_\varepsilon^d, N] \leq C(d)(\log 1/\varepsilon)^{-d/2}$$

*and, for $\varepsilon < \eta(d)$,*

$$d_{\text{Kol}}[\widetilde{B}_\varepsilon^d, N] \geq c(d)(\log 1/\varepsilon)^{-d/2}.$$

PROOF. We denote by

$$\widetilde{\kappa}_j(d,\varepsilon), \qquad j=1,2,\ldots,$$

the sequence of cumulants of the random variable $\widetilde{B}_\varepsilon^d$. We deal separately with the cases $d=1$ and $d \geq 2$.

(*Case* $d = 1$.) In this case, **W** is a standard Brownian motion on $[0,1]$ so that $\widetilde{B}^1_\varepsilon$ takes the form $\widetilde{B}^1_\varepsilon = I_2(f_\varepsilon)$, where $I_2$ is the double Wiener–Itô integral with respect to **W** and

$$(5.3) \qquad f_\varepsilon(x,y) = (4\log 1/\varepsilon)^{-1/2}[(x \vee y \vee \varepsilon)^{-1} - 1].$$

Lengthy (but standard) computations yield the following estimates: as $\varepsilon \to 0$,

$$\widetilde{\kappa}_2(1,\varepsilon) \longrightarrow 1,$$

$$\widetilde{\kappa}_j(1,\varepsilon) \asymp \left(\log \frac{1}{\varepsilon}\right)^{1-j/2}, \qquad j \geq 3.$$

The conclusion now follows from Proposition 3.8.

(*Case* $d \geq 2$.) In this case, $\widetilde{B}^d_\varepsilon$ has the form $\widetilde{B}^d_\varepsilon = I_2(f^d_\varepsilon)$, with

$$(5.4) \qquad f^d_\varepsilon(x_1,\ldots,x_d;y_1,\ldots,y_d) = (4\log 1/\varepsilon)^{-d/2} \prod_{j=1}^d [(x_j \vee y_j \vee \varepsilon)^{-1} - 1].$$

By using (3.13), one sees that

$$(2^{j-1}(j-1)!)^{-1} \times \widetilde{\kappa}_j(d,\varepsilon) = [(2^{j-1}(j-1)!)^{-1} \times \widetilde{\kappa}_j(1,\varepsilon)]^d,$$

so the conclusion follows once again from Proposition 3.8. $\square$

**6. Exact asymptotics in the Breuer–Major CLT.** Let $B$ be a fractional Brownian motion (fBm) with Hurst index $H \in (0, \frac{1}{2})$, that is, $\{B_x : x \geq 0\}$ is a centered Gaussian process with covariance given by

$$E(B_x B_y) = \tfrac{1}{2}(x^{2H} + y^{2H} - |x-y|^{2H}), \qquad x,y \geq 0.$$

It is well known that, for every choice of the parameter $H \in (0, \frac{1}{2})$, the Gaussian space generated by $B$ can be identified with an isonormal Gaussian process of the type $X = \{X(h) : h \in \mathfrak{H}\}$, where the real and separable Hilbert space $\mathfrak{H}$ is defined as follows: (i) denote by $\mathscr{E}$ the set of all $\mathbb{R}$-valued step functions on $\mathbb{R}_+$; (ii) define $\mathfrak{H}$ as the Hilbert space obtained by closing $\mathscr{E}$ with respect to the scalar product

$$\langle \mathbf{1}_{[0,x]}, \mathbf{1}_{[0,y]}\rangle_{\mathfrak{H}} = E(B_x B_y).$$

Such a construction implies, in particular, that $B_x = X(\mathbf{1}_{[0,x]})$. The reader is referred to, for example, [21] for more details on fBm, including crucial connections with fractional operators. We also define $\rho(\cdot)$ to be the covariance function associated with the stationary process $x \mapsto B_{x+1} - B_x$, that is,

$$\begin{aligned}\rho(x) &:= E[(B_{t+1} - B_t)(B_{t+x+1} - B_{t+x})] \\ &= \tfrac{1}{2}(|x+1|^{2H} + |x-1|^{2H} - 2|x|^{2H}), \qquad x \in \mathbb{R}.\end{aligned}$$

Now, fix an *even* integer $q \geq 2$ and set

$$Z_T := \frac{1}{\sigma(T)\sqrt{T}} \int_0^T H_q(B_{u+1} - B_u)\, du, \qquad T > 0,$$

where $H_q$ is the $q$th Hermite polynomial defined in (2.14) and where

$$\sigma(T) := \sqrt{\mathrm{Var}\left(\frac{1}{\sqrt{T}} \int_0^T H_q(B_{u+1} - B_u)\, du\right)} = \sqrt{\frac{q!}{T} \int_{[0,T]^2} \rho^q(u-v)\, du\, dv}.$$

Observe that each $Z_T$ can be represented as a multiple Wiener–Itô integral of order $q$ and also that

$$\sigma^2(T) \xrightarrow[T \to \infty]{} \sigma^2(\infty) := q! \int_{-\infty}^{+\infty} \rho^q(x)\, dx < +\infty.$$

According to, for example, the main results in [3] or [9], one always has the following CLT:

$$Z_T \xrightarrow[T \to \infty]{\mathrm{Law}} Z \sim \mathscr{N}(0,1)$$

(which also holds for odd values of $q$). The forthcoming Theorem 6.1 shows that the techniques of this paper may be used to deduce an exact asymptotic relation (as $T \to \infty$) for the difference $P(Z_T \leq z) - \Phi(z)$, where $\Phi(z) = P(N \leq z)$ $[N \sim \mathscr{N}(0,1)]$. We stress that the main results of this section deal with the case of a generic Hermite polynomial of even order $q \geq 2$, implying that our techniques even provide explicit results outside the framework of *quadratic* functionals, such as those analyzed in Sections 4 and 5. In what follows, we use the notation

$$\widehat{\sigma}^2(\infty) := \frac{q^2}{\sigma^4(\infty)} \sum_{s=1}^{q-1} (s-1)!^2 \binom{q-1}{s-1}^4 (2q-2s)!$$

(6.1)
$$\times \int_{\mathbb{R}^3} \rho^s(x_1)\rho^s(x_2)\rho^{q-s}(x_3)$$
$$\times \rho^{q-s}(x_2 + x_3 - x_1)\, dx_1\, dx_2\, dx_3$$

and

(6.2) $$\widehat{\gamma}(\infty) = -\frac{q!(q/2)!\binom{q}{q/2}^2}{2\sigma^3(\infty)} \int_{\mathbb{R}^2} \rho^{q/2}(x)\rho^{q/2}(y)\rho^{q/2}(x-y)\, dx\, dy.$$

THEOREM 6.1. *There exists a constant $C > 0$ such that*

(6.3) $$d_{\mathrm{Kol}}(Z_T, N) = \sup_{z \in \mathbb{R}} |P(Z_T \leq z) - \Phi(z)| \leq \frac{C}{\sqrt{T}}.$$

*Moreover, for any fixed $z \in \mathbb{R}$, we have*

(6.4) $$\sqrt{T}(P(Z_T \leq z) - \Phi(z)) \xrightarrow[T \to \infty]{} \frac{\widehat{\gamma}(\infty)}{3}(z^2 - 1)\frac{e^{-z^2/2}}{\sqrt{2\pi}}.$$

PROOF. The proof is divided into three steps.

STEP 1. Let us first prove the following convergence:

(6.5) $$\sqrt{T}\left(\frac{1}{q}\|DZ_T\|_{\mathfrak{H}}^2 - 1\right) \xrightarrow[T \to \infty]{\text{Law}} \mathcal{N}(0, \widehat{\sigma}^2(\infty)),$$

where $\widehat{\sigma}^2(\infty)$ is given by (6.1). Note that, once (6.5) is proved to be true, one deduces immediately that, as $T \to \infty$,

$$\text{Var}\left(\frac{1}{q}\|DZ_T\|_{\mathfrak{H}}^2 - 1\right) \sim \frac{\widehat{\sigma}^2(\infty)}{T},$$

so that (6.3) follows from Theorem 2.4. Now, to prove that (6.5) holds, start by using the well-known relation between Hermite polynomials and multiple integrals to write

$$H_q(B_{u+1} - B_u) = I_q(\mathbf{1}_{[u,u+1]}^{\otimes q}).$$

As a consequence, we have

$$DZ_T = \frac{q}{\sigma(T)\sqrt{T}} \int_0^T I_{q-1}(\mathbf{1}_{[u,u+1]}^{\otimes q-1})\mathbf{1}_{[u,u+1]}\,du.$$

Thus, by an appropriate use of the multiplication formula (2.4), one has that

$$\|DZ_T\|_{\mathfrak{H}}^2 = \frac{q^2}{\sigma^2(T)T}\int_{[0,T]^2}\rho(u-v)I_{q-1}(\mathbf{1}_{[u,u+1]}^{\otimes q-1})I_{q-1}(\mathbf{1}_{[v,v+1]}^{\otimes q-1})\,du\,dv$$

$$= \frac{q^2}{\sigma^2(T)T}\int_{[0,T]^2}\sum_{r=0}^{q-1}r!\binom{q-1}{r}^2 I_{2q-2-2r}(\mathbf{1}_{[u,u+1]}^{\otimes q-1-r} \otimes \mathbf{1}_{[v,v+1]}^{\otimes q-1-r})$$

$$\times \rho^{r+1}(u-v)\,du\,dv$$

$$= \frac{q^2}{\sigma^2(T)T}\sum_{s=1}^{q}(s-1)!\binom{q-1}{s-1}^2$$

$$\times \int_{[0,T]^2} I_{2q-2s}(\mathbf{1}_{[u,u+1]}^{\otimes q-s} \otimes \mathbf{1}_{[v,v+1]}^{\otimes q-s})\rho^s(u-v)\,du\,dv,$$

yielding

$$\frac{1}{q}\|DZ_T\|_{\mathfrak{H}}^2 - 1$$

$$= \frac{q}{\sigma^2(T)T} \sum_{s=1}^{q-1} (s-1)! \binom{q-1}{s-1}^2$$

$$\times \int_{[0,T]^2} I_{2q-2s}(\mathbf{1}_{[u,u+1]}^{\otimes q-s} \otimes \mathbf{1}_{[v,v+1]}^{\otimes q-s}) \rho^s(u-v) \, du \, dv.$$

We shall first prove that, for every $s \in \{1, \ldots, q-1\}$,

(6.6)
$$\frac{1}{\sqrt{T}} \int_{[0,T]^2} I_{2q-2s}(\mathbf{1}_{[u,u+1]}^{\otimes(q-s)} \otimes \mathbf{1}_{[v,v+1]}^{\otimes(q-s)}) \rho^s(u-v) \, du \, dv$$
$$\xrightarrow[T\to\infty]{\text{Law}} \mathcal{N}(0, \widehat{\sigma}_s^2(\infty)),$$

where

$$\widehat{\sigma}_s^2(\infty) := (2q-2s)! \int_{\mathbb{R}^3} \rho^s(x_1) \rho^s(x_2) \rho^{q-s}(x_3) \rho^{q-s}(x_2 + x_3 - x_1) \, dx_1 \, dx_2 \, dx_3.$$

Fix $s \in \{1, \ldots, q-1\}$. First, observe that

$$\widehat{\sigma}_s^2(T) := \mathrm{Var}\left( \frac{1}{\sqrt{T}} \int_{[0,T]^2} I_{2q-2s}(\mathbf{1}_{[u,u+1]}^{\otimes(q-s)} \otimes \mathbf{1}_{[v,v+1]}^{\otimes(q-s)}) \rho^s(u-v) \, du \, dv \right)$$

$$= \frac{(2q-2s)!}{T} \int_{[0,T]^4} \rho^s(u-v) \rho^s(w-z) \rho^{q-s}(u-w)$$
$$\times \rho^{q-s}(v-z) \, du \, dv \, dw \, dz$$

$$\xrightarrow[T\to\infty]{} \widehat{\sigma}_s^2(\infty)$$

so that (6.6) holds if and only if the following convergence takes place:

(6.7)
$$Q_T^{(s)} := \frac{1}{\widehat{\sigma}_s(T)\sqrt{T}} \int_{[0,T]^2} I_{2q-2s}(\mathbf{1}_{[u,u+1]}^{\otimes(q-s)} \otimes \mathbf{1}_{[v,v+1]}^{\otimes(q-s)}) \rho^s(u-v) \, du \, dv$$
$$\xrightarrow[T\to\infty]{\text{Law}} \mathcal{N}(0,1).$$

We have

$$DQ_T^{(s)} = \frac{2q-2s}{\widehat{\sigma}_s(T)\sqrt{T}}$$
$$\times \int_{[0,T]^2} \rho^s(u-v) I_{2q-2s-1}(\mathbf{1}_{[u,u+1]}^{\otimes(q-s-1)} \otimes \mathbf{1}_{[v,v+1]}^{\otimes(q-s)}) \mathbf{1}_{[u,u+1]} \, du \, dv.$$

Thus, $\|DQ_T^{(s)}\|_{\mathfrak{H}}^2$ is given by

$$\frac{(2q-2s)^2}{\widehat{\sigma}_s^2(T)T} \int_{[0,T]^4} \rho^s(u_1-u_2) \rho^s(u_3-u_4) \rho(u_1-u_3)$$

$$\times I_{2q-2s-1}(\mathbf{1}_{[u_1,u_1+1]}^{\otimes(q-s-1)} \otimes \mathbf{1}_{[u_2,u_2+1]}^{\otimes(q-s)})$$

$$\times I_{2q-2s-1}(\mathbf{1}_{[u_3,u_3+1]}^{\otimes(q-s-1)} \otimes \mathbf{1}_{[u_4,u_4+1]}^{\otimes(q-s)}) \, du_1 \cdots du_4$$

$$= \frac{(2q-2s)^2}{\widehat{\sigma}_s^2(T)T}$$

$$\times \int_{[0,T]^4} \rho^s(u_1-u_2)\rho^s(u_3-u_4)\rho(u_1-u_3)$$

$$\times \left( \sum_{t=0}^{q-s} t! \binom{2q-2s-1}{t}^2 \right.$$

$$\times I_{4q-4s-2-2t}(\mathbf{1}_{[u_1,u_1+1]}^{\otimes(q-s-1)} \otimes \mathbf{1}_{[u_3,u_3+1]}^{\otimes(q-s-1)}$$

$$\otimes \mathbf{1}_{[u_2,u_2+1]}^{\otimes(q-s-t)} \otimes \mathbf{1}_{[u_4,u_4+1]}^{\otimes(q-s-t)})$$

$$\left. \times \rho^t(u_2-u_4) \right) du_1 \cdots du_4$$

$$+ \frac{(2q-2s)^2}{\widehat{\sigma}_s^2(T)T}$$

$$\times \int_{[0,T]^4} \rho^s(u_1-u_2)\rho^s(u_3-u_4)\rho(u_1-u_3)$$

$$\times \left( \sum_{t=q-s+1}^{2q-2s-1} t! \binom{2q-2s-1}{t}^2 \right.$$

$$\times I_{4q-4s-2-2t}(\mathbf{1}_{[u_1,u_1+1]}^{\otimes(2q-2s-1-t)}$$

$$\otimes \mathbf{1}_{[u_3,u_3+1]}^{\otimes(2q-2s-1-t)})$$

$$\left. \times \rho^{q-s}(u_2-u_4)\rho^{t-q+s}(u_1-u_3) \right) du_1 \cdots du_4.$$

Consequently, $\frac{1}{2q-2s}\|DQ_T^{(s)}\|_{\mathfrak{H}}^2 - 1$ is given by

$$\frac{2q-2s}{\widehat{\sigma}_s^2(T)T}$$

$$\times \int_{[0,T]^4} \rho^s(u_1-u_2)\rho^s(u_3-u_4)\rho(u_1-u_3)$$

$$\times \left( \sum_{t=1}^{q-s+1} (t-1)! \binom{2q-2s-1}{t-1}^2 \right.$$

$$\times I_{4q-4s-2t}(\mathbf{1}_{[u_1,u_1+1]}^{\otimes(q-s-1)} \otimes \mathbf{1}_{[u_3,u_3+1]}^{\otimes(q-s-1)}$$

$$\otimes \mathbf{1}_{[u_2,u_2+1]}^{\otimes(q-s+1-t)} \otimes \mathbf{1}_{[u_4,u_4+1]}^{\otimes(q-s+1-t)})$$

$$\times \rho^{t+1}(u_2 - u_4) \Big) du_1 \cdots du_4$$

$$+ \frac{2q-2s}{\hat{\sigma}_s^2(T)T}$$

$$\times \int_{[0,T]^4} \rho^s(u_1 - u_2)\rho^s(u_3 - u_4)\rho(u_1 - u_3)$$

$$\times \left( \sum_{t=q-s+2}^{2q-2s-1} (t-1)! \binom{2q-2s-1}{t-1}^2 \right.$$

$$\times I_{4q-4s-2t}(\mathbf{1}_{[u_1,u_1+1]}^{\otimes(2q-2s-t)}$$

$$\otimes \mathbf{1}_{[u_3,u_3+1]}^{\otimes(2q-2s-t)})$$

$$\times \rho^{q-s}(u_2 - u_4)\rho^{t-q+s-1}(u_1 - u_3) \Big) du_1 \cdots du_4.$$

For a fixed $t$ such that $1 \le t \le q-s+1$, we have that

$$E \Big| \frac{1}{T} \int_{[0,T]^4} \rho^s(u_1 - u_2)\rho^s(u_3 - u_4)\rho(u_1 - u_3)\rho^{t+1}(u_2 - u_4)$$

$$\times I_{4q-4s-2t}(\mathbf{1}_{[u_1,u_1+1]}^{\otimes(q-s-1)} \otimes \mathbf{1}_{[u_3,u_3+1]}^{\otimes(q-s-1)}$$

$$\otimes \mathbf{1}_{[u_2,u_2+1]}^{\otimes(q-s+1-t)} \otimes \mathbf{1}_{[u_4,u_4+1]}^{\otimes(q-s+1-t)}) du_1 \cdots du_4 \Big|^2$$

$$= \frac{1}{T^2} \int_{[0,T]^8} \rho^s(u_1-u_2)\rho^s(u_3-u_4)\rho(u_1-u_3)$$

$$\times \rho^t(u_2-u_4)\rho^s(u_5-u_6)\rho^s(u_7-u_8)$$

$$\times \rho(u_5-u_7)\rho^t(u_6-u_8)\rho^{q-s-1}(u_1-u_5)$$

$$\times \rho^{q-s-1}(u_3-u_7)\rho^{q-s+1-t}(u_2-u_6)$$

$$\times \rho^{q-s+1-t}(u_4-u_8) du_1 \cdots du_8$$

$$\underset{T\to\infty}{\sim} \frac{1}{T} \int_{\mathbb{R}^7} \rho^s(x_1)\rho^s(x_2)\rho(x_3)\rho^t(x_2+x_3-x_1)\rho^s(x_4)\rho^s(x_5)\rho(x_6)$$

$$\times \rho^t(x_5+x_6-x_4)\rho^{q-s-1}(x_7)$$

$$\times \rho^{q-s-1}(x_6+x_7-x_3)\rho^{q-s+1-t}(x_4+x_7-x_1)$$

$$\times \rho^{q-s+1-t}(x_5 + x_6 + x_7 - x_2 - x_3)\, dx_1 \cdots dx_7$$

tends to zero as $T \to \infty$. Similarly, we can prove, for a fixed $t$ such that $q - s + 2 \leq t \leq 2q - 2s - 1$, that

$$E\left| \frac{1}{T} \int_{[0,T]^4} \rho^s(u_1 - u_2)\rho^s(u_3 - u_4)\rho^{t-q+s}(u_1 - u_3)\rho^{q-s}(u_2 - u_4) \right.$$
$$\left. \times I_{4q-4s-2t}(\mathbf{1}_{[u_1, u_1+1]}^{\otimes(2q-2s-t)} \otimes \mathbf{1}_{[u_3, u_3+1]}^{\otimes(2q-2s-t)})\, du_1 \cdots du_4 \right|^2$$

tends to zero as $T \to \infty$. Thanks to the main result in [22], the last relation implies that, for each $s$, the convergence (6.7) holds and, therefore, (6.6) is verified. Finally, by combining (6.6) with the results in [25] and [26], we obtain (6.5). Indeed, by using the orthogonality and isometric properties of multiple stochastic integrals, we can write

$$\mathrm{Var}\left(\sqrt{T}\left(\frac{1}{q}\|DZ_T\|_{\mathfrak{H}}^2 - 1\right)\right)$$
$$= \frac{q^2}{\sigma^4(T)T} \sum_{s=1}^{q-1} (s-1)!^2 \binom{q-1}{s-1}^4 (2q-2s)!$$
$$\times \left\langle \int_{[0,T]^2} \mathbf{1}_{[u,u+1]}^{\otimes(q-s)} \otimes \mathbf{1}_{[v,v+1]}^{\otimes(q-s)} \rho^s(u-v)\, du\, dv, \right.$$
$$\left. \int_{[0,T]^2} \mathbf{1}_{[w,w+1]}^{\otimes(q-s)} \otimes \mathbf{1}_{[z,z+1]}^{\otimes(q-s)} \rho^s(w-z)\, dw\, dz \right\rangle_{\mathfrak{H}^{\otimes(2q-2s)}}$$
$$= \frac{q^2}{\sigma^4(T)T} \sum_{s=1}^{q-1} (s-1)!^2 \binom{q-1}{s-1}^4 (2q-2s)!$$
$$\times \int_{[0,T]^4} \rho^s(u-v)\rho^s(w-z)$$
$$\times \rho^{q-s}(u-w)\rho^{q-s}(v-z)\, du\, dv\, dw\, dz,$$

from which we easily deduce that $\mathrm{Var}(\sqrt{T}(\frac{1}{q}\|DZ_T\|_{\mathfrak{H}}^2 - 1)) \xrightarrow[T \to \infty]{} \widehat{\sigma}^2(\infty)$.

STEP 2. Let us prove the following convergence:

(6.8) $$\left(Z_T, \sqrt{T}\left(\frac{1}{q}\|DZ_T\|_{\mathfrak{H}}^2 - 1\right)\right) \xrightarrow[T \to \infty]{\text{Law}} (U, V)$$

with $(U, V)$ a centered Gaussian vector such that $E(U^2) = 1$,

$$E(V^2) = \widehat{\sigma}^2(\infty)$$

and

$$E(UV) = -\widehat{\gamma}(\infty)$$
$$= \frac{q!(q/2)!\binom{q}{q/2}^2}{2\sigma^3(\infty)} \int_{\mathbb{R}^2} \rho^{q/2}(x)\rho^{q/2}(y)\rho^{q/2}(x-y)\, dx\, dy.$$

Observe that we already know that $Z_T \overset{\text{Law}}{\to} U$ and also that (6.5) is verified. Note, also, that we have proven (6.5) by first decomposing $\sqrt{T}(q^{-1}\|DZ_T\|_{\mathfrak{H}}^2 - 1)$ into a finite sum of multiple integrals and then by showing that each multiple integral satisfies an appropriate CLT. As a consequence, according to part B of Theorem 2.6 [with $G_n$ replaced by $\sqrt{T}(q^{-1}\|DZ_T\|_{\mathfrak{H}}^2 - 1)$], it is sufficient to show the following convergence:

(6.9)
$$E\bigg(Z_T \times \sqrt{T}\bigg(\frac{1}{q}\|DZ_T\|_{\mathfrak{H}}^2 - 1\bigg)\bigg)$$
$$\xrightarrow[T\to\infty]{} \frac{q!(q/2)!\binom{q}{q/2}^2}{2\sigma^3(\infty)} \int_{\mathbb{R}^2} \rho^{q/2}(x)\rho^{q/2}(y)\rho^{q/2}(x-y)\, dx\, dy.$$

By the orthogonality of multiple stochastic integrals, we can write

$$E\bigg(Z_T \times \sqrt{T}\bigg(\frac{1}{q}\|DZ_T\|_{\mathfrak{H}}^2 - 1\bigg)\bigg)$$
$$= \frac{q}{\sigma^3(T)T}\bigg(\frac{q}{2}-1\bigg)!\binom{q-1}{\frac{q}{2}-1}^2$$
$$\quad\times \int_{[0,T]^3} \rho^{q/2}(u-v) E(I_q(\mathbf{1}_{[w,w+1]}^{\otimes q}) I_q(\mathbf{1}_{[u,u+1]}^{\otimes q/2} \otimes \mathbf{1}_{[v,v+1]}^{\otimes q/2}))\, du\, dv\, dw$$
$$= \frac{qq!}{\sigma^3(T)T}\bigg(\frac{q}{2}-1\bigg)!\binom{q-1}{\frac{q}{2}-1}^2 \int_{[0,T]^3} \rho^{q/2}(u-v)\rho^{q/2}(u-w)$$
$$\qquad\qquad\qquad\qquad\qquad\qquad\qquad \times \rho^{q/2}(w-v)\, du\, dv\, dw$$
$$\xrightarrow[T\to\infty]{} \frac{qq!}{\sigma^3(\infty)}\bigg(\frac{q}{2}-1\bigg)!\binom{q-1}{\frac{q}{2}-1}^2 \int_{\mathbb{R}^2} \rho^{q/2}(x)\rho^{q/2}(y)\rho^{q/2}(x-y)\, dx\, dy$$
$$= \frac{q!(q/2)!\binom{q}{q/2}^2}{2\sigma^3(\infty)} \int_{\mathbb{R}^2} \rho^{q/2}(x)\rho^{q/2}(y)\rho^{q/2}(x-y)\, dx\, dy.$$

STEP 3. Step 1 and convergence (6.8) imply that, as $T \to \infty$,

$$\varphi(T) \sim \frac{\widehat{\sigma}(\infty)}{\sqrt{T}},$$

where $\varphi(T) = \text{Var}(1 - q^{-1}\|DZ_T\|_{\mathfrak{H}}^2 - 1)$ and

$$\left(Z_T, \frac{1 - 1/q\|DZ_T\|_{\mathfrak{H}}^2}{\varphi(T)}\right) \xrightarrow[T \to \infty]{\text{Law}} (U, \widehat{\sigma}(\infty)^{-1}V).$$

As a consequence, one can apply Theorem 3.1 in the case $\rho = \frac{\widehat{\gamma}(\infty)}{\widehat{\sigma}(\infty)}$ (the remaining assumptions are easily verified), yielding that

$$\varphi(T)^{-1}(P(Z_T \leq z) - \Phi(z)) \xrightarrow[T \to \infty]{} \frac{\widehat{\gamma}(\infty)}{3\widehat{\sigma}(\infty)}(z^2 - 1)\frac{e^{-z^2/2}}{\sqrt{2\pi}}$$

from which the conclusion follows. $\square$

**Acknowledgments.** We are grateful to D. Marinucci for discussions about Edgeworth expansions. We thank an anonymous referee for suggesting that one can prove Theorem 2.4 without assuming that $F$ has an absolutely continuous distribution, as well as for other insightful remarks.

Laboratoire de Probabilités
  et Modèles Aléatoires
Université Pierre et Marie Curie
Boîte courrier 188
4 Place Jussieu
75252 Paris Cedex 5
France
E-mail: ivan.nourdin@upmc.fr

Equipe Modal'X
Université Paris Ouest–Nanterre la Défense
200 Avenue de la République
92000 Nanterre
France
E-mail: giovanni.peccati@gmail.com